\documentclass[12pt,bezier]{amsart}
\usepackage{amsfonts,latexsym,amscd}
\newtheorem{theorem}{Theorem}[section]
\newtheorem{definition}[theorem]{Definition}
\newtheorem{proposition}[theorem]{Proposition}
\newtheorem{corollary}[theorem]{Corollary}
\newtheorem{lemma}[theorem]{Lemma}
\newtheorem{conjecture}[theorem]{Conjecture}
\newtheorem{vanishingconjecture}[theorem]{Vanishing Conjecture}

\newcommand{\thm}[1]{Theorem \ref{#1}}
\newcommand{\prop}[1]{Proposition \ref{#1}}
\newcommand{\lma}[1]{Lemma \ref{#1}}

\newcommand{\con}[1]{Conjecture \ref{#1}}
\newcommand{\vancon}[1]{Vanishing Conjecture \ref{#1}}

\newcommand{\rarrow}{\rightarrow}

\newcommand{\calA}{{\mathcal A}}

\newcommand{\calS}{{\mathcal S}}

\newcommand{\calV}{{\mathcal V}}
\newcommand{\calR}{{\mathcal R}}

\begin{document}

\title{On Simplicial Commutative Rings 
with Vanishing Andr\'e-Quillen Homology}
\author{James M. Turner}
\address{Department of Mathematics\\
College of the Holy Cross\\
One College Street\\
Worcester, MA 01610-2395}
\email{jmturner@math.holycross.edu}
\thanks{Research was partially supported by an NSF-NATO postdoctoral fellowship}
\date{August, 1997}



\begin{abstract}
We propose a generalization of a conjecture of D. Quillen, on the
vanishing of Andr\'e-Quillen homology, to simplicial commutative 
rings. This conjecture characterizes a notion of local complete 
intersection, extended to the simplicial setting, under a suitable 
hypothesis on the local characteristic. Further, under the condition 
of finite-type homology, we then prove the conjecture in the case of 
a simplicial commutative algebra augmented over a field of non-zero 
characteristic. As a consequence, we obtain a proof of Quillen's 
conjecture for a Noetherian commutative algebra - again augmented over a 
field of non-zero characteristic.
\end{abstract}

\maketitle

\section{Introduction}

In \cite{And} and \cite{Qui2}, M. Andr\'e and D. Quillen constructed the
notion of a homology\\ $D_*(A|R;M)$ for a commutative algebra $A$,
over a ring $R$, and an $A$-module $M$.  It was then conjectured
(see Section 5 of \cite{Qui2}) that, under suitable conditions on
$R$ and $A$, the vanishing of the homology in sufficiently
high degrees determines $A$ as a local complete intersection. In
particular, for local rings, the conjecture takes the following form.

\begin{conjecture}\label{con1.1}
Let $R$ be a (Noetherian, commutative) local ring with residue
field $\Bbb{F}$, and let $D_s(\Bbb{F}|R) = D_s(\Bbb{F}|R;\Bbb{F}), \ s
\geq 0$.  Then the following are equivalent:
\begin{itemize}
\item[(1)]
$D_s(\Bbb{F}|R) = 0, \quad s >> 0$;
\item[(2)] $D_s(\Bbb{F}|R) = 0, \quad s \geq 3$;
\item[(3)] $R$ is a complete intersection.
\end{itemize}
\end{conjecture}

In this form,  \ref{con1.1} was proven by L. Avramov, in outline form,
in \cite{Avr}, and, in much greater generality, in \cite{Avr2}.
Recall that a local ring $R$ is a complete intersection if its $I$-adic
completion $\hat{R}$ is a quotient of a complete regular ring by an
ideal generated by a regular sequence.  From this description and the
properties of Andr\'e-Quillen homology, the implications (3)
$\Rightarrow$ (2) $\Rightarrow$ (1) in \ref{con1.1} are immediate.

The objective of this paper is to extend a version of \ref{con1.1} for
simplicial local rings in an effort to bring the full power of
simplicial homotopy theory to bear on this type of problem and thereby
obtain a different proof of \ref{con1.1} closer in spirit to the
topological results of J.-P. Serre in \cite{Serre} and Y. Umeda in
\cite{Ume}.

In \cite{Qui3}, D. Quillen gave a construction of Andr\'e-Quillen
homology $D_*(A|B;M)$ where $B$ is a simplicial commutative ring, $A$ a
simplicial commutative $B$-algebra, and $M$ a simplicial $A$-module.

Let $\calR_{\Bbb{F}}$ be the category of (commutative) local rings, with
residue field $\Bbb{F}$, and $s \calR_{\Bbb{F}}$ the category of
simplicial objects over $\calR_{\Bbb{F}}$. It follows from \cite{Qui1}
that $s \calR_{\Bbb{F}}$ has three classes of maps, called weak
equivalences, fibrations, and cofibrations, giving it the structure of a
closed simplicial model category.  Using this structure, we say that a
simplicial local ring $R$ is an {\em  $n$-extension} if there
is a cofibration sequence
$$
\Sigma^{n-1} S_0 \rarrow R \rarrow S_{1},
$$
in the homotopy category $Ho(s \calR_{\Bbb{F}})$, such that $S_{0}$ 
is polynomial in $\calR_{\Bbb{F}}$ and $\hat{S_{1}} \cong 
\Sigma^{n}\bar{S}_{1}$ in $Ho(s \calR_{\Bbb{F}})$.  Here $\Sigma$ denotes
the suspension in $Ho(s \calR_{\Bbb{F}})$.  (See $\S$I.2 and $\S$I.3 of
\cite{Qui1} for the theory of suspension and cofibration sequences in
homotopical algebra.)

\begin{definition}\label{def1.2}
Let $R$ be an object of $s \calR_{\Bbb{F}}$. Then:
\begin{itemize}
\item[(1)] $R$ is {\bf regular} if $R$ is a 1-extension
with $\bar{S}_{1}$ smooth in $\calR_{\Bbb{F}}$.
\item[(2)] $R$ is a {\bf complete intersection} if it is a
1-extension with $\bar{S}_1$ regular in $s \calR_{\Bbb{F}}$.
\item[(3)] $R$ is {\bf $Q$-bounded} if $(Q \pi_* R)_s = 0$ for $s >> 0$,
and {\em bounded} if  $\pi_s R = 0$ for $s \gg 0$.
\item[(4)] The {\bf simplicial dimension} of $R$ is the integer
$$
s \cdot \dim R = \max\{s| \, D_s(\Bbb{F}|R) \neq 0\}.
$$
We then say that $R$ has {\bf finite simplicial dimension} if $s \cdot
\dim R < \infty$.
\item[(5)] $R$ is said to have {\bf finite-type homology} provided each
$D_q(\Bbb{F}|R)$ is a finite dimensional $\Bbb{F}$-vector space.
\item[(6)] If $R$ has both finite-type homology and finite simplicial
dimension, we call $R$ {\bf finite}.
\end{itemize}
Given a simplicial commutative ring $R$, then $R$ is said to be {\em
locally} of any one of (1) -- (6) provided $R_{\wp}$ is such, for each
simplicial prime ideal $\wp$ in $R$. Given a simplicial prime ideal $\wp$
in $R$, we denote by $\Bbb{F}(\wp)$ the residue field of $R_{\wp}$ and we say
that $R$ is {\em locally of non-zero characteristic} provided
char$\Bbb{F}(\wp)\neq 0$ for all such $\wp$.
\end{definition}

We can now state our proposed simplicial generalization of \con{con1.1}.

\begin{vanishingconjecture}\label{vancon1.3}
Let $R$ be a locally finite simplicial commutative ring which is locally
of non-zero characteristic. Then $R$ is a locally complete intersection
if and only if $R$ is locally $Q$-bounded.
\end{vanishingconjecture}

In the rational case, while a complete intersection may be both
$Q$-bounded and of finite simplicial dimension, the converse is not true.
See the note following \prop{prop3.7} regarding counter-examples.

To demonstrate the validity of \ref{vancon1.3}, we consider the
subcategory, $\calA_{\Bbb{F}}$, of $\calR_{\Bbb{F}}$, consisting of
augmented $\Bbb{F}$-algebras, i.e., unitary $\Bbb{F}$-algebras $A$
together with a fixed $\Bbb{F}$-algebra map $A \rarrow \Bbb{F}$, called
the augmentation of $A$.

In this paper, we give evidence for \vancon{vancon1.3} by proving:

\begin{theorem}\label{thm1.5}
Suppose $A$ is a finite simplicial augmented commutative
$\Bbb{F}$-algebra with char$\Bbb{F} > 0$. Then $A$ is bounded if and
only if $A$ is a complete intersection.
\end{theorem}

\begin{corollary}\label{cor1.5}
Let $A$ be an augmented Noetherian commutative $\Bbb{F}$-algebra, 
$char\Bbb{F} > 0$. Then $A$ is a complete intersection, as a local algebra, 
if and only if $A$ has finite simplicial dimension.
\end{corollary}

{\em Proof.}  By IV.55 of \cite{And}, $H^{Q}_{*}(A)$ is of finite-type. 
Thus, by \thm{thm1.5}, $A$ has finite simplicial dimension if and
only if $A$ is a complete intersection, as a simplicial algebra, if
and only if $A$ is a complete intersection, as a local algebra, by
the classical implication of (2) $\Rightarrow$ (3) in \con{con1.1}
(see Proposition 26 of \cite{And}). \hfill $\Box$
\bigskip

\subsection*{Organization of this paper}

In this section, we review the closed simplicial model category
structure for $s \calA_{\Bbb{F}}$ and the construction and properties of
homotopy and Andr\'e-Quillen homology. In Section 3, we describe
the notion Postnikov envelopes for objects
of $s \calA_{\Bbb{F}}$ and explore its properties. In Section 4, we
study the homotopy of n-extensions. Finally, in Section 5, we
introduce and study the notion of a Poincar\'e series for a simplicial
algebra, obtaining just enough information to prove \thm{thm1.5}.

\subsection*{Acknowledgements}
The author would like to thank Haynes Miller, for suggesting this
project along with the direction it should take, and Jean Lannes for
many useful directions as well as for making his stay in France
worthwhile.  Most of the work on this project was done while the
author was visiting the Institut des Hautes \'Etudes Scientifique and
the Ecole Polytechnique. He would like to thank them for their
hospitality and use of their facilities during his stay. Finally, the
author would like to thank Julie Riddleburger for putting this paper
into LaTeX form.

\section{Homotopy Theory of Simplicial Augmented\\ Commutative Algebras}
\setcounter{equation}{0}

We now review the closed simplicial model category structure for $s
\calA_{\Bbb{F}}$.  We will assume the reader is familiar with the
general theory of homotopical algebra given in \cite{Qui1}.

We call a map $f: \ A \rarrow B$ in $s \calA_{\Bbb{F}}$ a
\begin{itemize}
\item[(i)] weak equivalence ($\stackrel{\sim}{\rightarrow}$)
$\Leftarrow\!\!\Rightarrow \pi_* f$ is an isomorphism;
\item[(ii)] fibration ($\rightarrow\!\!\!\!\rightarrow$)
$\Leftarrow\!\!\Rightarrow f$
surjects in positive degrees;
\item[(iii)] cofibration($\hookrightarrow$) $\Leftarrow \!\!\Rightarrow f$
is a retract of
an almost free map.
\end{itemize}

Here a map $f: \ A \rarrow B$ in $s \calA_{\Bbb{F}}$ is {\em almost
free} if there is an almost simplicial $\Bbb{F}$-vector space (no $d_0$)
$V$ (see \cite{Goe1}) together with a map of almost simplicial
$\Bbb{F}$-vector spaces $V \rarrow IB$ such that the induced map $A
\otimes S(V) \stackrel{\cong}{\longrightarrow} B$ is an isomorphism of
almost simplicial algebras. Here $S$ is the symmetric algebra functor.

Now given a finite simplicial set $K$ and a simplicial algebra $A$,
define $A \wedge K$ and $A^K$ by
$$
(A \wedge K)_n = \bigotimes_{K_{n}} A_n
$$
and
$$
(A^K)_n = \prod_{K_{n}} A_n.
$$
Here the tensor product $\otimes $ is the coproduct in $s
\calA_{\Bbb{F}}$.  The product in $s \calA_{\Bbb{F}}$ is defined as
$\Lambda \times_{\Bbb{F}} \Gamma$, for $\Lambda, \Gamma$ in $s
\calA_{\Bbb{F}}$, so that the diagram
$$
\begin{array}{ccc}
\Lambda \times_{\Bbb{F}}\Gamma & \longrightarrow & \Gamma \\[2mm]
\downarrow & & \hspace*{10pt} \downarrow \epsilon \\[2mm]
\Lambda & \longrightarrow & \Bbb{F} \\[-3mm]
        & \epsilon &
\end{array}
$$
is a pullback of simplicial vector spaces.

\begin{theorem}\label{thm2.1} (\cite{Qui1}, \cite{Mil}, and \cite{Goe1})
 With these definitions, $s \calA_{\Bbb{F}}$ is a closed simplicial
model category.
\end{theorem}

Given a simplicial vector space $V$, define its normalized chain complex
$NV$ by
\begin{equation}
N_nV = V_n/(\mbox{Im} s_0 + \cdots + \mbox{Im} s_n)
\end{equation}
and $\partial: \ N_n V \rarrow N_{n-1}V$ is $\partial = \sum^n_{i=0}
(-1)^id_i$. The homotopy groups $\pi_* V$ of $V$ is defined as
$$
\pi_n V = H_n(NV), \quad n \geq 0.
$$
Thus for $A$ in $s \calA_{\Bbb{F}}$ we define $\pi_* A$ as above.
The Eilenberg-Zilber theorem (see \cite{Mac}) shows that the algebra
structure on $A$ induces an algebra structure on $\pi_*A$.

If we let $\calV$ be the category of $\Bbb{F}$-vector spaces, then there
is an adjoint pair
$$
S: \ \calV \Leftarrow\!\!\Rightarrow \calA_{\Bbb{F}}: \, I,
$$
where $I$ is the augmentation ideal function and $S$ is the symmetric
algebra functor.  For an object $V$ in $\calV$ and $n \geq 0$, let
$K(V,n)$ be the associated Eilenberg-MacLane object in $s \calV$ so that
$$
\pi_s K(V,n) = \left\{\begin{array}{ll}
V & s = n; \\[2mm]
0 & s \neq n.
\end{array}\right.
$$
Let $S(V,n) = S(K(V,n))$, which is an object of $s \calA_{\Bbb{F}}$.

For $A$ in $\calA_{\Bbb{F}}$, the indecomposable functor $QA =
I(A)/I^2(A)$ which is an object of $\calV$.  Furthermore, we have an
adjoint pair
$$
Q: \ \calA_{\Bbb{F}} \Leftarrow \!\!\Rightarrow \calV: \ (-)_+
$$
where $V_+$, for $V$ in $\calV$, is the simplicial algebra $V \oplus
\Bbb{F}$ where
$$
(v,r) \cdot (w,s) = (sv+rw,rs)
$$
for $(v,r), (w,s) \in V \oplus \Bbb{F}$. $(-)_+, Q$ provides an
equivalence between $\calV$ and the category of abelian group objects
in $\calA_{\Bbb{F}}$.

For $A$ in $s \calA_{\Bbb{F}}$, we define its Andr\'e-Quillen homology,
as per \cite{Goe1} and \cite{Goe2}, by
$$
H^Q_s(A) = \pi_s QX, \quad s \geq 0,
$$
where we choose a factorization
$$
\Bbb{F} \hookrightarrow X \stackrel{\sim\hspace*{5pt}}{\rightarrow \!\!\!\!\!
\rightarrow} A
$$
of the unit $\Bbb{F} \rarrow A$ as a cofibration and a trivial
fibration.  This definition is independent of the choice of
factorization as any two are homotopic over $A$ (note that every object
of $s \calA_{\Bbb{F}}$ is fibrant). It is known (see, for example,
\cite{Mil}) that
$$
H^Q_s(A) = D_s(A|\Bbb{F};\Bbb{F}).
$$
From the transitivity sequence, one can easily check that
$D_0(\Bbb{F}|A) = 0$, and $D_{s+1}(\Bbb{F}|A)\\ \cong H^Q_s(A)$ for all
$s \geq 0$.

Now, as shown in \cite{Goe2},
$$
\pi_n A = [S(n),A],
$$
where $S(n) = S(\Bbb{F},n)$ and $[\quad,\quad]$ denotes the morphisms in
$Ho(s \calA_{\Bbb{F}})$.  Thus the primary operational structure for
the homotopy groups in $s \calA_{\Bbb{F}}$ is determined by $\pi_*
S(V_0)$ for any $V_0$ in $s \calV$.  By Dold's theorem \cite{Dold} there
is a triple $\calS$ on graded vector spaces so that
\begin{equation}
\pi_* S(V) \cong \calS(\pi_* V)
\end{equation}
encoding this structure. If char$\Bbb{F} = 0$, $\calS$ is the free skew
symmetric functor and, if char$\Bbb{F} > 0$, $\calS$ is the free divided
power algebra on the underlying vector space of a certain free algebra
constructed from the input (see, for example, \cite{Bou} and \cite{Goe1}).

Now recall that maps $A \stackrel{f}{\rightarrow} B
\stackrel{g}{\rightarrow} C$ is a {\em cofibration sequence} in $Ho(s
\calA_{\Bbb{F}})$ if $f$ is isomorphic to a cofibration $X
\stackrel{u}{\rightarrow} Y$, of cofibrant objects, with cofibre $Y
\stackrel{v}{\rightarrow} Z$ isomorphic to $g$. Thus given any map $f: \
A \rarrow B$ in $s \calA_{\Bbb{F}}$ there is a cofibration sequence $A
\stackrel{f}{\rarrow} B \rarrow M(f)$ in $Ho(s \calA_{\Bbb{F}})$ formed
by factoring $\Bbb{F} \rarrow A$ into $\Bbb{F} \hookrightarrow \bar{A}
\stackrel{\sim}{\rightarrow\!\!\!\!\rightarrow}  A$, form the diagram
$$
\begin{array}{ccc}
\bar{A} & \hookrightarrow & X \\[2mm]
s
\begin{picture}(1,1)\put(1,1) {$\downarrow$} \put(1,5)
{$\downarrow$}\end{picture}
\hspace*{10pt} &&
\begin{picture}(1,1)\put(1,1) {$\downarrow$} \put(1,5)
{$\downarrow$}\end{picture}
\hspace*{10pt} \wr \\[2mm]
A & \stackrel{f}{\rightarrow} & B
\end{array}
$$
and then let $M(f) = X \otimes_{\bar{A}}\Bbb{F}$, which is cofibrant.  As 
an example, the {\em suspension} $\Sigma A$ of an object
$A$ in $s \calA_{\Bbb{F}}$ by $M(\epsilon)$, where $\epsilon: \ A \rarrow
\Bbb{F}$ is the augmentation.

Finally, recall that the {\it completion} $\hat{A}$ of a simplicial augmented 
algebra $A$ is defined as
$$
\hat{A} = lim_{t} A/I^{t}.
$$
If $f: A \to B$ is a map of simplicial algebras, we denote by 
$\hat{f}: \hat{A} \to \hat{B}$ the induced map of completions.

We can now summarize methods for computing homotopy and Andr\'e-Quillen
homology that we will need for this paper.

\begin{proposition}\label{prop2.4}
\begin{itemize}
\item[(1)] If $f: \ A \stackrel{\sim}{\rarrow} B$ is a weak equivalence
in $s \calA_{\Bbb{F}}$, then $H^Q_*(f): \ H^Q_*(A)
\stackrel{\cong}{\rarrow} H^Q_*(B)$ is an isomorphism.
\item[(2)] for any $A$ in $s \calA_{\Bbb{F}}$ there is a spectral sequence
$$
E^1_{s,t} = \calS_{s} (H^Q_*(A)) \Rightarrow \pi_t \hat{A},
$$
called Quillen's Fundamental spectral sequence (see \cite{Qui2} and
\cite{Qui3}), which converges when $H^{Q}_{0}(A)=0$.
\item[(3)] There is a Hurewicz homomorphism $h: \ I\pi_* A \rarrow
H^Q_*(A)$ such that if $A$ is connected and $H^Q_s (A) = 0$ for $s<n$
then $A$ is $n$-connected and $h: \ \pi_n A \stackrel{\cong}{\rarrow}
H^Q_n(A)$ is an isomorphism.
\end{itemize}
\end{proposition}

{\em Proof.}
(1) is a standard result. See, for example, \cite{Qui2} or \cite{Goe2}.
For (2), see chapter IV of \cite{Goe1} and \cite{Tur}. Finally, (3) is 
in \cite{Goe1}. \hfill $\Box$
\bigskip

The following is a selection of results from \cite{Tur}.

\begin{proposition}\label{compprop}
Let $A$ and $B$ be in $s \calA_{\Bbb{F}}$. Then
\begin{itemize}
\item[(1)] if $f:A \to B$ is an $H^{Q}_{*}$-isomorphism then 
$\hat{f}: \hat{A} \to \hat{B}$ is a weak equivalence,
\item[(2)] if $A$ is connected then $\pi_{*}\hat{A} \cong \pi_{*}A$, 
\item[(3)] $H^{Q}_{*}(\hat{A}) \cong H^{Q}_{*}(A)$, and
\item[(4)] $Q\pi_{*}\hat{A} \cong Q\pi_{*}A$.
\end{itemize}
\end{proposition}

{\em Remark.} If $H^{Q}_{0}(A) = 0$ then \prop{compprop} (4) follows 
from a Quillen fundamental spectral sequence argument. This is 
due to the fact that while this spectral sequence doesn't directly 
converge to $\pi_{*}A$ it does allow, under the above condition on $H^{Q}_{0}$,  
sufficient information to be extracted about the indecomposables (see \cite{Tur} 
for further details). This case is sufficient for our needs.
\bigskip

\begin{proposition}\label{prop2.5}
Let $A \stackrel{f}{\rarrow} B \stackrel{g}{\rarrow} C$ be a cofibration
sequence in $Ho(s \calA_{\Bbb{F}})$.  Then:
\begin{itemize}
\item[(1)]
There is a long exact sequence
$$
\begin{array}{l}
\cdots \rarrow H^Q_{s+1} (C) \stackrel{\partial}{\rarrow} H^Q_s(A)
\stackrel{H^Q_*(f)}{\rarrow} H^Q_s(B) \\[3mm]
\hspace*{20pt} \stackrel{H^Q_*(g)}{\rarrow} H^Q_s(C)
\stackrel{\partial}{\rarrow} H^Q_{s-1}(C) \rarrow \cdots
\end{array}
$$
\item[(2)] There is a first quadrant spectral sequence of algebras
$$
E^2_{s,t} = Tor^{\pi_*A}_{s} (\pi_* B, \Bbb{F})_t
\Rightarrow \pi_{s+t} C,
$$
which we refer to as the Eilenberg-Moore spectral sequence.
\item[(3)] If $A$ is connected, there is a first quadrant spectral
sequence of algebras
$$
E^2_{s,t} = \pi_s (C \otimes \pi_t A) \Rightarrow \pi_{s+t} B,
$$
which we refer to as the Serre spectral sequence.
\item[(4)] If $A$ is connected and $C$ is $n$-connected, then there is a
homomorphism $\tau: \ \pi_{n+1} C \rarrow \pi_n A$, called the {\em
transgression}, such that the diagram
$$
\begin{array}{ccccccc}
\pi_{n+1}B & \stackrel{\pi_* f}{\rarrow} & \pi_{n+1}C & \stackrel{\tau}{\rarrow}
& \pi_nA & \stackrel{\pi_* f}{\rarrow} & \pi_n B \\[1mm]
h \downarrow \hspace*{8pt}
&& h \downarrow \hspace*{8pt}
&& \hspace*{8pt} \downarrow h && \hspace*{8pt} \downarrow h \\[1mm]
H^Q_{n+1}B & \stackrel{H^Q_*f}{\rarrow} & H^Q_{n+1}(C) &
\stackrel{\partial}{\rarrow} & H^Q_n(A) & \stackrel{H^Q_*f}{\rarrow} &
H^Q_n B
\end{array}
$$
commutes and the top sequence is exact.
\end{itemize}
\end{proposition}

{\em Proof.}
(1) is just the transitivity sequence for $H^Q_*$. See \cite{Goe1}.

(2) is the spectral sequence of Theorem 6(b) in $\S$II.6 of \cite{Qui1}.
See also \cite{Goe1}. By Theorem 6(d) in $\S$II.6 of \cite{Qui1}, there
is a $1^{st}$-quadrant spectral sequence
$$
E^2_{*,*} = \pi_* (B \otimes_A \pi_* A) \Rightarrow \pi_* B,
$$
where $\pi_*A$ is an $A$-module via the augmentation $A \rarrow \pi_0
A$.  Here we can assume our cofibration sequence is a cofibration with
cofibre $C$.  Since $A$ is connected, then $B \otimes_A \pi_* A \cong C
\otimes \pi_* A$.  The algebra structure follows from the construction
of the spectral sequence and the fact that $A \stackrel{f}{\rarrow} B$
is a map of simplicial algebras.  This gives us (3).

For (4), since $A$ is connected and $C$ is $n$-connected, then in the
Serre spectral sequence
$$
d^{n+1}: \ \pi_{n+1} C \cong E^{n+1}_{n+1,0} \rarrow E^{n+1}_{0,n} \cong
\pi_n A,
$$
which we propose is our desired map $\tau$.  From this same spectral
sequence, we have $\pi_s A \cong \pi_s B$, $s < n$, and, using methods
modified from the next section, we can assume that $N_s IC = 0$ for $s
\leq n$ and $N_{n+1}B \rarrow\!\!\!\!\rarrow N_{n+1}C$ is surjective.
Since $E^{1}_{n+1,0} = N_{n+1}C$ and $E^2_{n+1,0} \cong E^{n+1}_{n+1,0}
\cong \pi_{n+1}C$, then $d^{n+1}$ is constructable in precisely the same
way as the boundary map in homological algebra. Since we can assume
cofibrancy of our objects under consideration, then the diagram
$$
\begin{array}{ccccccc}
H_{n+1}(NB) & \rarrow & H_{n+1}(NC) & \stackrel{d^{n+1}}{\rarrow}
& H_n(NA) & \rarrow & H_n(NB) \\
\downarrow && \downarrow && \downarrow && \downarrow \\
H_{n+1}(NQB) & \rarrow & H_{n+1}(NQC) & \stackrel{\partial}{\rarrow} &
H_n(NQA) & \rarrow & H_n(NQB)
\end{array}
$$
commutes by naturality. The result follows. \hfill $\Box$

\section{Postnikov Envelopes}
\setcounter{equation}{0}

In this section, we construct and determine some
properties of a useful tool for studying simplicial algebras.
First, we recall the following standard result which will be useful for
us (see section II.4 of \cite{Qui1}).

\begin{lemma}\label{lma3.4}
Let $V$ and $W$ be simplicial vector spaces.  Then the map
$$
[V,W] \rarrow \mbox{Hom}_{\calV_{*}}(\pi_* V, \pi_* W)
$$
is an isomorphism.
\end{lemma}

\begin{proposition}\label{prop3.5}
Let $A$ in $s \calA_{\Bbb{F}}$. Then
\begin{itemize}
\item[(1)] There is a 
map of simplicial algebras 
$$
f_{0}:S(H^{Q}_{0}(A),0)\to A
$$ 
which induces an isomorphism on $H^{Q}_{0}$.
\item[(2)] Suppose $A$ is $(n-1)$-connected for $n \geq 1$.
Then there exists a map in $s \calA_{\Bbb{F}}$,
$$
f_n: \ S(H^Q_nA,n) \rarrow A,
$$
which is an isomorphism on $\pi_n$ and $H^Q_n$.
\end{itemize}
\end{proposition}

{\em Proof.} 
(1) Let $\iota: H^{Q}_{0}(A)\to I\pi_{0}A$ be a choice of splitting for 
the surjection $I\pi_{0}A \to H^{Q}_{0}(A)$. By \lma{lma3.4}, $f$ can be 
chosen to be the adjoint of 
the map of simplicial vector spaces $K(H^{Q}_{0}(A),0)\to IA$ induced 
by $\iota$. By the transitivity sequence, $H^{Q}_{0}(A(1))=0$ so the 
fundamental spectral sequence for $A(1)$ converges, by \prop{prop2.4} 
(2), so $A(1)$ is connected.
(2) By the Hurewicz theorem, \prop{prop2.4} (3), the map
$h: \ \pi_nA \rarrow H^Q_n A$ is an isomorphism.  Now the adjoint functors
$$
S: \ s \calV \Leftarrow\!\!\Rightarrow s \calA_{\Bbb{F}}: \ I
$$
induce an adjoint pair
$$
S: \ Ho(s \calV) \Leftarrow\!\!\Rightarrow Ho(s \calA_{\Bbb{F}}): \ I.
$$
Thus we have isomorphisms
\begin{eqnarray*}
[S(H^Q_nA,n),A]
& \cong & [K(H^Q_nA,n),IA] \\[2mm]
& \cong & \mbox{Hom}_{\calV}(H^Q_n A, \pi_n IA),
\end{eqnarray*}
using \lma{lma3.4}. Choosing $f_n$ to correspond to the inverse of $h$
gives the result. \hfill $\Box$
\bigskip

Now given $A$ we form the {\em Postnikov envelopes} as the
sequence of cofibrations
$$
A(1) \stackrel{j_{2}}{\hookrightarrow} A(2)
\stackrel{j_{3}}{\hookrightarrow} \cdots
\stackrel{j_{n}}{\hookrightarrow}
A(n) \stackrel{j_{n+1}}{\hookrightarrow} \cdots
$$
with the following properties: 
\begin{itemize}
\item[(1)] $A(1) = \widehat{M(f_{0})}$,
\item[(2)] for each $n\geq 1$, $A(n)$ is a $(n-1)$-connected and for $s \geq n$,
$$
H^Q_s A(n) \cong H^Q_sA.
$$
\item[(3)] There is a cofibration sequence
$$
S(H^Q_{n} A, n) \rarrow A(n) \stackrel{j_{n+1}}{\rarrow} A(n+1).
$$
\end{itemize}
The existence of a Postnikov envelopes follows easily from 
\prop{compprop}, \prop{prop3.5}, and

\begin{lemma}\label{lma3.6}
If $A$ is $(n-1)$-connected, for $n\geq 1$, then the cofibre $M(f_n)$ of $f_n: \
S(H^Q_nA,n) \rarrow A$ is $n$-connected and satisfies $H^Q_sM(f_n) \cong
H^Q_sA$ for $s > n$.
\end{lemma}

{\em Proof.}
This follows from \ref{prop3.5} and the transitivity sequence 
$$
H^Q_{s+1}M(f_n) \rarrow H^Q_s S(H^Q_n A,n) \rarrow H^Q_s A \rarrow
H^Q_s M(f_n).
$$
\hfill $\Box$
\bigskip

{\bf Note:}  We have been implicitly using the computation
$$
H^Q_s S(V,n) = \pi_s QS(V,n) = \pi_s K(V,n) = V
$$
for $s = n$ and 0 otherwise.  The converse holds as well.

\begin{proposition}\label{prop3.7}
Let $A$ be connected in $s \calA_{\Bbb{F}}$ and suppose $H^Q_s A = 0, \
s \neq n > 0$.  Then $A \cong S(H^Q_n A,n)$ in $Ho(s \calA_{\Bbb{F}})$.
\end{proposition}

{\em Proof.}  Since $A$ is connected, then $A$ is $(n-1)$-connected by
the Hurewicz theorem. By \ref{prop3.5}, $f_n: \ S(H^Q_nA,n) \rarrow A$ is
an $H^Q_n$-isomorphism and hence a weak equivalence by
\ref{prop2.4}(2). \hfill $\Box$
\bigskip

{\bf Note.} From this proposition, if char$\Bbb{F}$=0 then $S(V,n)$
has simplicial dimension $n$ and $\pi_{*}S(V,n)$ is free
skew-commutative on a basis of $V$ concentrated in degree $n$. Thus
$S(V,n)$ is $Q$-bounded, for any $n$, showing that \vancon{vancon1.3} fails
in the zero characteristic case.
\bigskip

\section{The Homotopy and Homology of $n$-Extensions}
\setcounter{equation}{0}

Call an object $A$ in $s \calA_{\Bbb{F}}$ a {\em simple
$n$-extension} if $A$ is an $n$-extension in $s
\calA_{\Bbb{F}}$ with $\bar{S}_1 = S(V_1,0)$, $V_1$ in $\calV$. 
Also, for this section and the next, we define the {\em simplicial
dimension} of $A$ to be
$$
s \cdot \dim A = \max\{s| \, H^{Q}_{s}(A) \neq 0\}
$$

We now proceed to prove:


\begin{theorem}\label{thm1.4}
Let $A$ be in $s \calA_{\Bbb{F}}$.  Then:
\begin{itemize}
\item[(1)] If $A$ is a connected simple $n$-extension for $n \geq 2$,
then, in $Ho(s \calA_{\Bbb{F}})$, we have
$$
A \cong S(H^Q_{n-1}(A),n-1) \otimes S(H^Q_n(A),n).
$$
\item[(2)] $A$ is a complete intersection if and only if $A$ is a simple
1-extension.
\item[(3)] If $A$ is a complete intersection then $H^Q_s(A) = 0$
for $s \geq 2$ and if $H^{Q}_{*}(A)$ is of finite-type then $A$ is
bounded.
\item[(4)] If $H^{Q}_{0}(A)=0$ and $H^{Q}_{*}(A)$ is of finite-type
then $\pi_{*}(\hat{A})$ is of finite-type.
\item[(5)] The Postnikov envelope $A(1)$ has the following properties:
\begin{itemize}
\item[(a)] If $A$ has finite simplicial dimension, then so does $A(1)$;
\item[(b)] If $H^{Q}_{0}(A)$ is finite and $\pi_{*}A$ is bounded 
then $\pi_{*}A(1)$ is $Q$-bounded.
\item[(c)] If $H^{Q}_{*}(A)$ is of finite-type then $H^{Q}_{*}(A(1))$ 
is also of finite-type.
\end{itemize}
\end{itemize}
\end{theorem}

We begin with

\begin{lemma}\label{lma4.1}
Let $A$ in $s \calA_{\Bbb{F}}$ be a connected simple $n$-extension for $n
\geq 2$.  Then $A$ is an $n$-extension
of the form
$$
S(H^Q_{n-1} A,n-1) \rarrow A \rarrow S(H^Q_n A,n).
$$
\end{lemma}

{\em Proof.}  Let $V_0, V_1$ be vector spaces so that there is a
cofibration sequence
$$
S(V_0, n-1) \rarrow A \rarrow S(V,n).
$$
Then the transitivity sequence tells us that $H^Q_s A = 0$, $s \neq n, \
n-1$ and there is an exact sequence
$$
0 \rarrow H^Q_nA \rarrow V_1 \rarrow V_0 \rarrow H^Q_{n-1}A \rarrow 0.
$$
Thus $A$ is $n-2$ connected and Postnikov tower gives us a cofibration sequence
$$
S(H^Q_{n-1}A, n-1) \rarrow A \rarrow A(n-1) = S(H^Q_n A,n).
\eqno\Box
$$
\bigskip

{\em Proof of \thm{thm1.4} (1).} By \lma{lma4.1}, there is a cofibration
sequence
$$
S(H^Q_{n-1},A,n-1) \stackrel{i}{\rarrow} A \stackrel{j}{\rarrow} S(H^Q_nA,n),
$$
where we can assume $A$ is cofibrant, $i$ is a cofibration, and $j$ is
the cofibre.  Consider the commuting diagram
$$
\begin{array}{ccc}
[S(H^Q_nA,n),A] & \stackrel{j_*}{\longrightarrow} &
[S(H^Q_nA,n),S(H^Q_nA,n)] \\[1mm]
\cong \downarrow \hspace*{10pt}
&&
\hspace*{10pt} \downarrow \cong \\[1mm]
[K(H^Q_n A,n),IA] && [K(H^Q_n A,n), IS(H^Q_n A,n)] \\[1mm]
\cong \downarrow \hspace*{10pt}
&&
\hspace*{10pt} \downarrow \cong \\[1mm]
\mbox{Hom}(H^Q_nA,I \pi_n A) & \stackrel{h_*}{\longrightarrow}
& \mbox{Hom}(H^Q_n A, H^Q_n A)
\end{array}
$$
Then $j$ will split, up to homotopy, if we can show that $h: \ \pi_n A
\rarrow H^Q_n A$ is onto.

By \prop{prop2.5}(4), there is a commutative diagram
$$
\begin{array}{ccccccc}
\pi_n A & \stackrel{\pi_* j}{\rarrow} & \pi_n S(H^Q_n A,n)
& \stackrel{\tau}{\rarrow} & \pi_{n-1} S(H^Q_{n-1},A, n-1) &
\stackrel{\pi_* i}{\rarrow} & \pi_{n-1} A \\[2mm]
h \downarrow \hspace*{10pt} && \cong \downarrow \hspace*{10pt} &&
\hspace*{10pt} \downarrow \cong && \hspace*{10pt} \downarrow \cong \\[2mm]
H^Q_n A & \stackrel{\cong}{\rarrow} & H^Q_n A & \stackrel{\partial = 0}{\rarrow}
& H^Q_{n-1}A & \stackrel{\cong}{\rarrow} & H^Q_{n-1}A
\end{array}
$$
with the rows exact. Thus $\pi_n j$ is onto and, hence, $h: \ \pi_n A
\rarrow H^Q_n A$ is onto. \hfill $\Box$
\bigskip

\begin{lemma}\label{lma4.3}
Suppose $A$ in $s \calA_{\Bbb{F}}$ is regular. Then
$S(H^{Q}_{0}(A),0) \cong \hat{A}$ in $Ho(s \calA_{\Bbb{F}})$.
\end{lemma}

{\em Proof.} By the standard transitivity sequence for $D_*$ applied to
$\Bbb{F} \rarrow A \rarrow \Bbb{F}$, $D_0(\Bbb{F}|A) = 0$ and
$D_{s+1}(\Bbb{F}|A) \cong H^Q_s(A)$, so since $A$ is regular, then
$H^Q_s(A) = 0$, $s > 0$. Thus $f_{0}$ is an $H^{Q}_{*}$-isomorphism 
and so $\hat{f_{0}}$ is a weak equivalence by \prop{compprop}. \hfill $\Box$
\bigskip

{\em Proof of \thm{thm1.4} (2).} If $A$ is a complete intersection
then it is a 1-extension of the form
$$
S_{0}\rightarrow A \rightarrow  S_{1}
$$
with $S_{0}$ polynomial and $\bar{S}_{1}$ regular as simplicial augmented 
algebra. By \lma{lma4.3}, A is thus a simple 1-extension. The converse is
clear. \hfill $\Box$
\bigskip

{\em Proof of \thm{thm1.4} (3).} If $A$ is a complete intersection,
then $H^Q_s(A) = 0, \ s \geq 2$ follows (2) and the transitivity
sequence. Consider now the Eilenberg-Moore spectral sequence
$$
E^{2}_{s,t}=Tor^{S(H^{Q}_{0}(A))}_{s}(\pi_{*}A,\Bbb{F})_{t}\Longrightarrow
\pi_{s+t}M(f_{0})
$$
which is a first quadrant homology-type spectral sequence of algebras.
Since $H^{Q}_{0}(A)$ is finite, $S(H^{Q}_{0}(A))$ has finite flat 
dimension and since, by \prop{compprop},
$Q\pi_{*}M(f_{0}) = Q\pi_{*}A(1) = Q\pi_{*}S(H^{Q}_{1}(A),1)$ is 
finite concentrated in degree 1 then
$\pi_{*}M(f_{0})$ is bounded and we can conclude, by an induction on 
$dim_{\Bbb{F}}H^{Q}_{0}(A)$, that $\pi_{*}A$ is bounded.
\hfill $\Box$
\bigskip

{\em Example.} Suppose an augmented commutative $\Bbb{F}$-algebra $B$ is a
complete intersection.  Then there is a complete regular algebra $\Gamma$ and 
an ideal
$I$, generated by a regular sequence, so that $\Gamma/I \cong \hat{B}$.  As we
saw, $\Gamma \cong S(V_0)$,  so the condition of regularity on $I$ is
equivalent to there being
a {\em projective extension}, that is, (see \cite{Goe1}) an extension
$$
\Bbb{F} \rarrow S(V_1) \stackrel{i}{\rarrow} S(V_0) \rarrow \hat{B} \rarrow 
\Bbb{F},
$$
so that $i$ makes $S(V_0)$ into a projective $S(V_1)$-module.  In
$Ho(s \calA_{\Bbb{F}})$, $\hat{B}$ is equivalent $M(i)$ and so there is
a cofibration sequence of the form
$$
S(V_0,0) \rarrow \hat{B} \rarrow S(V_1,1).
$$
Thus, $\hat{B}$, and hence $B$, is a complete intersection as a simplicial 
algebra.
\hfill $\Box$
\bigskip

{\em Proof of \thm{thm1.4} (4).} Since $H^{Q}_{0}(A)=0$, the
fundamental spectral sequence
$$
E^1_{s,t} = \calS_{s} (H^Q_*(A))_{t} \Rightarrow \pi_t \hat{A},
$$
converges. From the known structure of $\calS$ (see e.g. \cite{Bou}),
if $V$ is a finite-dimensional vector space then each
$\calS_{s}(V)_{t}$ is finite and $\calS_{s}(V)_{t}=0, s \gg 0$ for
each fixed t. The result follows. \hfill $\Box$
\bigskip

%

{\em Proof of \thm{thm1.4} (5).} First, (a) is immediate from the
transitivity sequence. For (b), $V = H^{Q}_{0}(A)$ is finite and the
Eilenberg-Moore spectral sequence has the form
$$
E^{2}_{s,t}=Tor^{S(V)}_{s}(\pi_{*}A,\Bbb{F})_{t}\Longrightarrow
\pi_{s+t}M(f_{0})
$$
Since $S(V)$ has finite flat dimension and $\pi_{*}A$ is a graded
$S(V)$-module then
$$
Tor^{S(V)}_{s}(\pi_{*}A,\Bbb{F})_{t} = Tor^{S(V)}_{s}(\pi_{t}A,\Bbb{F})
$$
vanishes for $s \gg 0$ and vanishes for $t \gg 0$ if $\pi_{*}A$ is
bounded. We conclude $\pi_{*}M(f_{0})$ is bounded and hence 
$\pi_{*}A(1)$ is $Q$-bounded, by \prop{compprop}.  Finally, for (c),
\lma{lma3.6} tells us that $H^{Q}_{s}(A)=H^{Q}_{s}(A(1))$ for
$s \geq 1$. Thus if  $H^{Q}_{*}(A)$ is of finite-type
then $H^{Q}_{*}(A(1))$ is of finite-type. \hfill $\Box$
\bigskip

\section{The Poincar\'e Series of a Simplicial Algebra}
\setcounter{equation}{0}

For this section, we assume char$\Bbb{F} = p > 0$.

Let $A$ be a connected simplicial augmented commutative
$\Bbb{F}$-algebra such that $\pi_{*}A$ is of finite-type. We define its
{\em Poincar\'e series} by
$$
\vartheta(A,t) = \sum_{n\geq 0}(dim_{\Bbb{F}}\pi_{n}A)t^{n}.
$$
If $V$ is a finite-dimensional vector space and $n>0$ we write
$$
\vartheta(V,n,t) = \vartheta(S(V,n),t).
$$
Given power series $f(t) = \sum a_{i}t^{i}$ and $g(t) = \sum b_{i}t^{i}$
we define the relation $f(t) \leq g(t)$ provided $a_{i}\leq
b_{i}$ for each $i\geq 0$.

\begin{lemma}\label{poiprop}
Given a cofibration sequence
$$
A \rightarrow B \rightarrow C
$$
of connected objects in $\calA_{\Bbb{F}}$ with finite-type homotopy
groups, then
$$
\vartheta(B,t) \leq \vartheta(A,t)\vartheta(C,t)
$$
which is an equality if the sequence is split.
\end{lemma}

{\em Proof.} From the Serre spectral sequence
$$
E^{2}_{s,t}=\pi_{s}(C\otimes\pi_{t}A) \Longrightarrow \pi_{s+t}B
$$
we have
$$
\vartheta(A,t)\vartheta(C,t) =
\sum_{n}(\sum_{s+t=n}dim_{\Bbb{F}}E^{2}_{s,t})t^{n} \geq
\vartheta(B,t).
$$
If the cofibration sequence is split then the spectral sequence
collapses, giving an equality. \hfill $\Box$
\bigskip

If $\Pi$ is a finitely-generated abelian group and $n>0$ let
$$
\vartheta(\Pi,n,t) = \sum_{s}(dim_{\Bbb{F}}H_{s}(K(\Pi,n);\Bbb{F}))t^{s}.
$$

\begin{lemma}\label{lmapoi}
Let $V$ be a finite-dimensional vector space and $\Pi$ a free abelian
group of the same dimension. Then for any $n>0$
$$
\vartheta(V,n,t) = \vartheta(\Pi,n,t).
$$
\end{lemma}

{\em Proof.} As shown in \cite{Car}, there is a weak equivalence of
simplicial vector spaces
$$
S(V,n) \rightarrow \Bbb{F}[K(\Pi,n)]
$$
which gives us the desired result. \hfill $\Box$
\bigskip

\begin{proposition}\label{phiprop}
Given a finite-dimensional vector space $V$ and any $n>0$ the
Poin-\\car\'e series $\vartheta(V,n,t)$ converges in the open unit disc.
\end{proposition}

{\em Proof.} This follows from \lma{lmapoi} and the results of J.P.
Serre in \cite{Serre} and Y. Umeda in \cite{Ume}. \hfill $\Box$
\bigskip

Now given two power series $f(t)$ and $g(t)$ we say $f(t) \sim g(t)$
provided $lim_{t\to \infty }f(t)/g(t)\\ = 1$. Given a Poincar\'e series
$\vartheta(V,n,t)$, for a finite-dimensional $\Bbb{F}$-vector space $V$
and $n>0$, let
$$
\varphi(V,n,t) = log_{p}\vartheta(V,n,1-p^{-t}).
$$

\begin{proposition}\label{poieq}
For $V$ an $\Bbb{F}$-vector space of finite dimension
$q$ and $n>0$ then $\varphi(V,n,t)$ converges on the real line and
$$
\varphi(V,n,t)\sim qt^{n-1}/(n-1)!.
$$
\end{proposition}

{\em Proof.} This follows from \lma{lmapoi} and Th\'eor\`eme 9b in
\cite{Serre}, for char$\Bbb{F}=2$,  and its generalization in \cite{Ume}.
\hfill $\Box$
\bigskip

A major step in proving \thm{thm1.5} will be accomplished with

\begin{theorem}\label{connprop}
Let $A$ be a connected finite simplicial augmented commutative
$\Bbb{F}$-algebra. Then if $A$ is $Q$-bounded we have $A \cong
S(H^{Q}_{1}(A),1)$ in $Ho(s \calA_{\Bbb{F}})$.
\end{theorem}

{\em Proof.} By \thm{thm1.4} (4) and \prop{compprop}, $\pi_{*}A$ is of 
finite-type and
hence, as it is also $Q$-bounded, bounded as well. Let $n = s \cdot \dim
A$. We must show that $n = 1$.

Consider the Postnikov envelope
$$
S(H^{Q}_{s-1}(A),s-1) \rightarrow A(s-1) \rightarrow
A(s)
$$
for each s. From the theory of cofibration sequences (see section I.3 of
\cite{Qui1}) the above sequence extends to a cofibration sequence
$$
A(s-1) \rightarrow A(s) \rightarrow
S(H^{Q}_{s-1}(A),s).
$$
Thus, by \lma{poiprop}, we have
$$
\vartheta(A(s),t) \leq
\vartheta(A(s-1),t)\vartheta(H^{Q}_{s-1}(A),s,t).
$$
Starting at $s=n-1$ and iterating this relation, we arrive at the
inequality
$$
\vartheta(A(n-1),t) \leq \vartheta(A,t)\prod_{s=1}^{n-2}
\vartheta(H^{Q}_{s}(A),s+1).
$$
Now, $A(n) \cong S(H^{Q}_{n}(A),n)$ by \prop{prop3.7},
but, by \thm{thm1.4} (1) and \lma{poiprop}, we have
$$
\vartheta(A(n-1),t) =
\vartheta(H^{Q}_{n-1}(A),n-1,t)\vartheta(H^{Q}_{n}(A),n,t).
$$
Since $\pi_{*}(A)$ is of finite-type and
bounded then there exists a $D>p$ such that
$\vartheta(A,t) \leq D$, in the open unit disc. Combining, we have
$$
\vartheta(H^{Q}_{n-1}(A),n-1,t)\vartheta(H^{Q}_{n}(A),n,t)
\leq D\prod_{s=1}^{n-2}\vartheta(H^{Q}_{s}(A),s+1).
$$
Applying a change of variables and $log_{p}$ to the above inequality, we get
$$
\varphi(H^{Q}_{n-1}(A),n-1,t)+\varphi(H^{Q}_{n}(A),n,t)
\leq d+\sum_{s=1}^{n-2}\varphi(H^{Q}_{s}(A),s+1).
$$
By \prop{poieq}, there is a polynomial $f(t)$ of degree $n-2$,
non-negative integer $a$, and positive integers $b$ and $d$ such that
$$
at^{n-2}+bt^{n-1} \leq d+f(t)
$$
which is clearly false for $n>1$. Thus $n=1$. The rest of the proof
follows from \prop{prop3.7}. \hfill $\Box$
\bigskip

{\em Proof of \thm{thm1.5}.} The ``if'' part is \thm{thm1.4} (3). We
thus concentrate on the ``only if'' part. We are given a finite simplicial
augmented commutative $\Bbb{F}$-algebra $A$ such that $\pi_{*}A$ is 
bounded. By \thm{thm1.4} (4) and (5), $A(1)$ is connected, $Q$-bounded, 
finite, and $\pi_{*}A(1)$ is of finite-type. Thus $A(1)$ is bounded as 
well and we conclude $A(1) \cong S(H^{Q}_{1}(A), 1)$ by \thm{connprop}. 
\hfill $\Box$

\end{document}